\RequirePackage{fix-cm}
\RequirePackage{amsmath}
\documentclass[smallextended]{svjour3}       

\usepackage{latexsym}
\usepackage{amsfonts}
\usepackage{graphicx} 
\usepackage{color}    
\usepackage{float}
\usepackage{longtable} 
\usepackage{adjustbox}
\usepackage{enumerate}
\usepackage{mathtools}
\usepackage{float}
\usepackage{multirow}
\usepackage{blkarray}
\usepackage{tikz}
\usepackage{pgfplots}
\pgfplotsset{compat=1.7}
\usepackage{pgfplotstable}
\usepackage{algorithm}
\usepackage[algo2e]{algorithm2e}
\usepackage{algpseudocode}

\newtheorem{thm}[lemma]{Theorem}

\newcommand{\rmv}[1]{}

\newcommand{\comment}[1]{}

\newcommand{\ftwo}{{\mathbb F}_{2}}
\newcommand{\ftwom}{{\mathbb F}_{2^m}}

\begin{document}

\title{A new class of irreducible pentanomials for polynomial based
multipliers in binary fields}

\titlerunning{A new class of irreducible pentanomials}


\author{Gustavo Banegas \and
        Ricardo Cust\'{o}dio \and Daniel Panario 
}

\institute{   Gustavo Banegas (\email{gustavo@cryptme.in}) at
              Technische Universiteit Eindhoven, Eindhoven, The Netherlands
           \and
              Ricardo Cust\'{o}dio (\email{ricardo.custodio@ufsc.br}) at
              Universidade Federal de Santa Catarina, Florian\'{o}polis, Brazil
           \and
              Daniel Panario (\email{daniel@math.carleton.ca}) at
              Carleton University, Ottawa, Canada. \\
}


\maketitle

\begin{abstract}
We introduce a new class of irreducible pentanomials over $\ftwo$ of the form
$f(x) = x^{2b+c} + x^{b+c} + x^b + x^c + 1$. Let $m=2b+c$ and use $f$ to
define the finite field extension of degree $m$. We give the exact number of
operations required for computing the reduction modulo $f$. We also provide
a multiplier based on Karatsuba algorithm in $\mathbb{F}_2[x]$
combined with our reduction process. We give the total cost of the multiplier and found that the
bit-parallel multiplier defined by this new class of polynomials
has improved XOR and AND complexity. Our multiplier has comparable time delay when compared to
other multipliers based on Karatsuba algorithm.
\end{abstract}

\section{Introduction}

Finite field extensions $\ftwom$ of the binary field $\ftwo$ play a central role in many engineering applications
and areas such as cryptography. Elements in these extensions are commonly represented using polynomial or normal bases.
We center in this paper on polynomial bases for bit-parallel multipliers.

When using polynomial bases, since $\ftwom \cong \ftwo[x]/(f)$
for an irreducible polynomial $f$ over $\ftwo$ of degree $m$,
we write elements in $\ftwom$ as polynomials over $\ftwo$ of
degree smaller than $m$. When multiplying with elements in $\ftwom$, a polynomial of degree up to $2m-2$ may arise.
In this case, a modular reduction is necessary to bring the resulting element back to $\ftwom$. Mathematically, any
irreducible polynomial can be used to define the extension.
In practice, however, the choice of the irreducible $f$ is crucial for fast and efficient field multiplication.

There are two types of multipliers in $\mathbb{F}_{2^m}$:
one-step algorithms and two-step algorithms. Algorithms of
the first type perform modular reduction while the elements
are being multiplied. In this paper, we are interested in
two-step algorithms, that is, in the first step the
multiplication of the elements is performed, and in the
second step the modular reduction is executed. Many
algorithms have been proposed for both types. An interesting
application of two-step algorithms is in several
cryptographic implementations that use the lazy reduction method \cite{unterluggauer2014efficient,bernstein2012neon}. For
example, in \cite{negre2013impact} it is shown the impact of lazy reduction
in operations for binary elliptic curves. An important application of the second part of our algorithm,
the reduction process, is to side-channel attacks. Indeed,
we prove that our modular reduction requires a constant
number of arithmetic operations, and as a consequence, it
prevents side-channel attacks.

The complexity of hardware circuits for finite field arithmetic in $\ftwom$ is related to the amount of space
and the time delay needed to perform the operations. Normally, the number of exclusive-or (XOR) and AND gates is
a good estimation of the space complexity. The time complexity is the delay due to the use of these gates.

Several special types of irreducible polynomials have been considered before, including polynomials with few nonzero
terms like trinomials and pentanomials (three and five nonzero terms, respectively), equally spaced
polynomials, all-one polynomials~\cite{Halbutogullari,itoh1989structure,Masoleh:2004}, and other special
families of polynomials~\cite{Wu2008}. In general, trinomials are preferred, but for
degrees where there are no irreducible trinomials, pentanomials are considered.

The analysis of the complexity using trinomials is known~\cite{wu2002bit}. However, there is no general
complexity analysis of a generic pentanomial in the literature. Previous results (see~\cite{FanHasan} for
 details) have focus on special classes of pentanomials, including:

\begin{itemize}

 \item $x^m + x^{b+1} + x^{b} + x^{b-1} + 1$, where $2 \leq b \leq m/2-1$ \cite{imana2013low,RodriguezKoc,zhang2001systematic,park2014new,imana2016high};

 \item $x^m + x^{b+1} + x^{b} + x + 1$, where $1 < b < m-1$ \cite{imana2006low,Masoleh:2004,RodriguezKoc,zhang2001systematic,park2014new,imana2016high};

 \item $x^{m} + x^{m-c} + x^{b} + x^{c} + 1$, where $1 \leq c < b < m-c$ \cite{Cilardo};

 \item $x^m + x^a + x^b + x^c + 1$, where $1 \leq c < b < a \leq m/2$ \cite{Masoleh:2004};

 \item $x^{m} + x^{m-s} + x^{m-2s} + x^{m-3s} + 1$, where $(m-1)/8 \leq s \leq (m-1)/3$ \cite{Masoleh:2004};

 \item $x^{4c} + x^{3c} + x^{2c} + x^{c} + 1$, where $c = 5^i$ and $i \geq 0$ \cite{Halbutogullari,Hasan}.

\end{itemize}

Like our family, these previous families focus on bit operations, i.e., operations that use only
 AND and XOR gates. In the literature it is possible to find studies that use computer words to perform the operations \cite{scott2007optimal,oliveira2016software} but this is not the focus of our work.

\subsection{Contributions of this paper}

In this paper, we introduce a new class of irreducible pentanomials with the following format:
\begin{equation} \label{ourpenta}
f(x) = x^{2b+c} + x^{b+c} + x^b + x^c + 1, b > c > 0.
\end{equation}
We compare our pentanomial with the first two families from the
list above. The reason to choose these two family is that
\cite{park2014new} presents a multiplier considering these families
with complexity 25\% smaller than the other existing works in the literature using quadratic algorithms. Since our multiplier is based on
Karatsuba's algorithm, we also compare our method with Karatsuba type algorithms.

An important reference for previously used polynomials and their
complexities is the recent survey on bit-parallel multipliers by
Fan and Hasan \cite{FanHasan}. Moreover, we observe that all
finite fields results used in this paper can be found in the
classical textbook by Lidl and Niederreiter \cite{LN}; see
\cite{HFF} for recent research in finite fields.

We prove that the complexity of the reduction depends on the exponents $b$ and $c$ of the pentanomial. A consequence of
 our result is that for a given degree $m=2b+c$, for any positive integers $b>c>0$, all irreducible polynomials in
our family have the same space and time complexity. We provide the exact number of XORs and gate delay required for
the reduction of a polynomial of degree $2m-2$ by our pentanomials. The number of XORs needed is $3m-2=6b+3c-2$ when $b \neq 2c$; for
$b=2c$ this number is $\frac{12}{5}m - 1=12c-1$. We also show that AND gates are not required in the reduction process.
It is easy to verify that our reduction algorithm is ``constant-time''
since it runs the same amount of operations independent of the inputs and
it avoids timing side-channel attacks \cite{fan2010state}.

For comparison purposes with other pentanomials proposed in the
literature, since the operation considered in those papers is the
product of elements in $\ftwom$, we also consider the number of ANDs
and XORs used in the multiplication prior to the reduction. In the
literature, one can find works that use the standard product or use some more efficient method of multiplication,
such as Karatsuba, and then add the complexity of the reduction.

In this paper, we use a Karatsuba multiplier combined with our fast reduction method. The total cost is
then $C m^{\log_2{3}} + 3m-2$ or $C m^{\log_2{3}} + \frac{12}{5}m-1$, depending on $b \not = 2c$ or $b = 2c$,
respectively. The constant $C$ of the Karatsuba multiplier depends on the implementation. In our experiments, $C$
is strictly less than $6$ for all practical degrees, up to degrees $1024$. For the reduction,
we give algorithms that achieve the above number of operations using any irreducible pentanomial in our family.
We compare the complexity of the Karatsuba multiplier with our reduction with the method proposed by Park et. al\cite{park2014new},
as well as, with Karatsuba variants given in~\cite{FanHasan}.

\subsection{Structure of the paper}
The structure of this paper is as follows. In Section~\ref{numberreductions}
we give the number of required reduction steps when using a pentanomial
$f$ from our family. We show that
for our pentanomials this number is $2$ or $3$. This fact is crucial
since such a low number of required reduction steps of our family allows
for not only an exact count of the XOR operations but also for a reduced time delay. Our strategy for that consists in describing the reduction process throughout equations, cleaning the redundant operations and presenting the final optimized algorithm. Section~\ref{reductions} provides
the first component of our strategy. In this section, we simply reduce a
polynomial of degree at most or exactly $2m-2$ to a polynomial of degree
smaller than $m$. The second component of our strategy is more delicate
and it allows us to derive the exact number of operations involved when
our pentanomial $f$ is used to define $\ftwom$.
Sections~\ref{twored}~and~\ref{threered} provide those analyses for the
cases when two and three steps of reduction are needed, that is, when
$c=1$ and $c>1$, respectively. We give algorithms and exact estimates
for the space and time complexities in those cases. Also, we describe a Karatsuba multiplier implementation combined
with our reduction.  In Section~\ref{sec:comparisons}, based on our implementation, we show that the number of XOR and AND gates
is better than the known space complexity in the literature. On the other hand, the time
complexity (delay) in our implementation is worse than quadratic methods but comparable with Karatsuba implementations. Hence, our multiplier
would be preferable in situations where space complexity and saving energy are more relevant than time complexity.
We demonstrate that our family contains many polynomials, including degrees where pentanomials are suggested by NIST. Conclusions
are given in Section~\ref{conclusions}.

\section{The number of required reductions}
\label{numberreductions}

When operating with two elements in $\ftwom$, represented by polynomials,
we obtain a polynomial of degree at most $2m-2$. In order to obtain the
corresponding element in $\ftwom$, a further division with remainder by
an irreducible polynomial $f$ of degree $m$ is required. We can see this
reduction as a process to bring the coefficient in interval $[2m-2,m]$ to a
position less than $m$. This  is done in steps. In each step, the
coefficients in interval $[2m-2,m]$ of the polynomial is substituted by the
equivalent bits following the congruence $x^m \equiv x^a+x^b+x^c+1$.
Once the coefficient in position $2m-2$ is brought to a position less
than $m$, the reduction is completed.

In this section, we carefully look into the number of steps needed to
reduce the polynomial by our polynomial $f$ given in Equation
(\ref{ourpenta}). The most important result of this section is that we
need at most $3$ steps of this reduction process using our polynomials.
This information is used in the next sections to give the exact number
of operations when the irreducible pentanomial given in Equation
(\ref{ourpenta}) is employed. This computation was possible because
our family has a small number of required reduction steps.

Let $ D_0(x) = \sum_{i=0}^{2m-2}d_ix^i$
be a polynomial over $\ftwo$. We want to compute $D_{red}$, the remainder of the division of $D_0$ by $f$, where $f$ has the form $f(x) = x^{2b+c} + x^{b+c} + x^b + x^c + 1$ with $2b+c=m$ and $b>c>0$.
The maximum number $k_a$ of reduction steps for a pentanomial
$x^{m} + x^{a} + x^b + x^c + 1$ in terms of the exponent
$a$ is given by Sunar and Ko\c{c} \cite{SunarKoc}
$$ k_a = \left \lfloor \frac{m-2}{m-a} \right \rfloor + 1.$$
In our case $m=2b+c$ and $a=b+c$, thus
\begin{equation} \label{eq:numero:reducoes}
  k_{b+c}  = \left\lfloor \frac{2b+c-2}{2b+c-b-c} \right\rfloor + 1 = \left\lfloor \frac{c-2}{b} \right\rfloor + 3 = \begin{cases}
        2 & \text{if } c = 1, \\
        3 & \text{if } c > 1.
     \end{cases}
\end{equation}
Using the same method as in \cite{SunarKoc}, we can derive the
number of steps required associated to the exponents
$b$ and $c$. These numbers are needed in Section~\ref{reductions}. We get
\begin{equation} \label{eq:numero:redb}
   k_{b}
   = \left\lfloor \frac{2b+c-2}{2b+c-b} \right\rfloor + 1
   = \left\lfloor \frac{b-2}{b+c} \right\rfloor + 2
   = 2,
\end{equation}
and
\begin{equation} \label{eq:numero:redc}
   k_c = \left\lfloor \frac{2b+c-2}{2b+c-c} \right\rfloor + 1
   = \left\lfloor \frac{c-2}{2b} \right\rfloor + 2 = \begin{cases}
        1 & \text{if } c = 1, \\
        2 & \text{if } c > 1.
     \end{cases}
\end{equation}
Thus, the reduction process for our family of pentanomials involves at
most three steps. This is a special property that our family enjoys.

The general process for the reduction proposed in this paper is given in the next section.
The special case $c=1$, that is when our polynomials have the form $f(x) = x^{2b+1} + x^{b+1} + x^b + x + 1$, requires two steps.
This family is treated in detail in Section~\ref{twored}.
The general case of our family $f(x) = x^{2b+c} + x^{b+c} + x^b + x^c + 1$ for $c>1$
involves three steps and is treated in Section~\ref{threered}.

\section{The general reduction process} \label{reductions}
\begin{figure}[htb]
   \centering
   \includegraphics[width = .4\columnwidth]{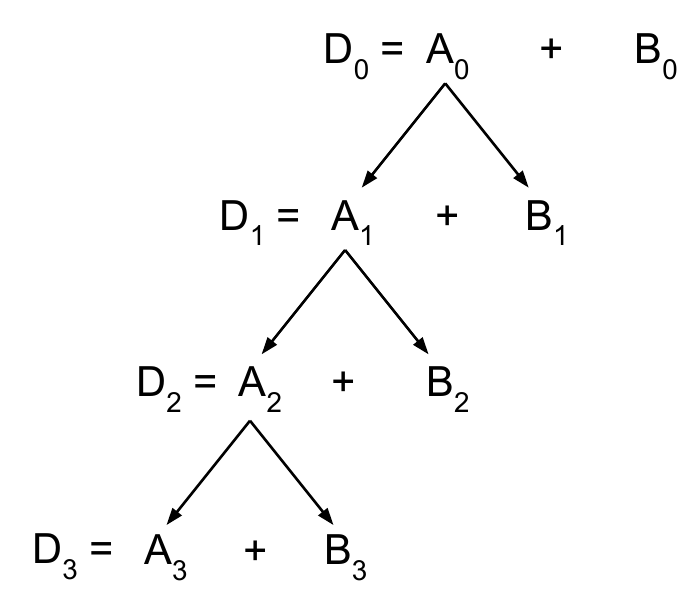}
   \caption{Tree representing the general reduction strategy.}
   \label{fig:arvore}
\end{figure}

The general process that we follow to get the original polynomial $D_0$
reduced to a polynomial of degree smaller than $m$ is depicted in
Figure~\ref{fig:arvore}. Without loss of generality, we consider the
polynomial to be reduced as always having degree $2m-2$. Indeed, the
cost to determine the degree of the polynomial to be reduced is
equivalent to checking if the leading coefficient is zero.

The polynomial $D_0$ to be reduced is split into two parts: $A_0$ is
the piece of the original polynomial with degree at least $m$ and hence
that requires extra work, while $B_0$ is formed by the terms of $D_0$
with exponents smaller than $m$ and so that it does not require to be
reduced. Dividing the leading term of $A_0$ by $f$ with remainder we
obtain $D_1$. In the same way as before, we split $D_1$ in two parts
$A_1$ and $B_1$ and repeat the process obtaining the tree of
Figure~\ref{fig:arvore}.

\subsection{Determining $A_0$ and $B_0$}
We trivially have
$$ D_0 (x) = A_0(x) + B_0(x)
   = \sum_{i=m}^{2m-2}d_ix^i + \sum_{i=0}^{m-1}d_ix^i,$$
and hence
\begin{equation} \label{eq:B0}
  A_0 = \sum_{i=m}^{2m-2}d_ix^i
        \quad \mbox{and } \quad
  B_0 = \sum_{i=0}^{m-1}d_ix^i.
\end{equation}

\subsection{Determining $A_1$ and $B_1$}
Using for clarity the generic form of a pentanomial over $\ftwo$,
$f(x) = x^{m} + x^{a} + x^b + x^c + 1$, dividing the leading term
of $A_0$ by $f$ and taking the remainder, we get
\begin{equation*}
D_1 = \sum_{i=0}^{m-2}d_{i+m}x^{i+a} + \sum_{i=0}^{m-2}d_{i+m}x^{i+b} +  \sum_{i=0}^{m-2}d_{i+m}x^{i+c} + \sum_{i=0}^{m-2}d_{i+m}x^i.
\end{equation*}

Separating the already reduced part of $D_1$ from the piece
of $D_1$ that still requires more work, we obtain
\begin{equation} \label{eq:A1}
A_1 = \sum_{i=m}^{m+a-2}d_{i+(m-a)}x^i + \sum_{i=m}^{m+b-2}d_{i+(m-b)}x^i + \sum_{i=m}^{m+c-2}d_{i+(m-c)}x^i,
\end{equation}
and
\begin{equation*}
B_1 = \sum_{i=a}^{m-1}d_{i+(m-a)}x^i + \sum_{i=b}^{m-1}d_{i+(m-b)}x^i +  \sum_{i=c}^{m-1}d_{i+(m-c)}x^i + \sum_{i=0}^{m-2}d_{i+m}x^i.
\end{equation*}
Since $m=2b+c$ and $a=b+c$, we have
\begin{equation*}
A_1 = \sum_{i=2b+c}^{3b+2c-2}d_{i+b}x^i + \sum_{i=2b+c}^{3b+c-2}d_{i+b+c}x^i +  \sum_{i=2b+c}^{2b+2c-2}d_{i+2b}x^i,
\end{equation*}
\begin{equation} \label{eq:B1}
B_1 = \sum_{i=b+c}^{2b+c-1}d_{i+b}x^i + \sum_{i=b}^{2b+c-1}d_{i+b+c}x^i + \sum_{i=c}^{2b+c-1}d_{i+2b}x^i  + \sum_{i=0}^{2b+c-2}d_{i+2b+c}x^i.
\end{equation}

\subsection{Determining $A_2$ and $B_2$}
As before, we divide the leading term of $A_1$ by $f$ and we obtain the remainder $D_2$. We get $D_2 = D_{2_a} + D_{2_b} + D_{2_c}$, where $D_{2_a}$, $D_{2_b}$ and $D_{2_c}$ refer to the reductions of the sums  in Equation~(\ref{eq:A1}).

We start with $D_{2_{a}}$:
$$ D_{2_{a}} = \sum_{i=0}^{a-2}d_{i+2m-a}x^i(x^a+x^b+x^c+1).$$
Separating $D_{2_a}$ in the pieces $A_{2_{a}}$ and
$B_{2_{a}}$, we get
$A_{2_{a}} = \sum_{i=m}^{2a-2}d_{i+2m-2a}x^i$ since $b+a-2 < m$,
and
\begin{equation*}
B_{2_{a}} = \sum_{i=a}^{m-1}d_{i+2m-2a}x^i + \sum_{i=b}^{a+b-2}d_{i+2m-a-b}x^i + \sum_{i=c}^{a+c-2}d_{i+2m-a-c}x^i + \sum_{i=0}^{a-2}d_{i+2m-a}x^i.
\end{equation*}
Substituting $m=2b+c$ and $a=b+c$, we get
$A_{2_{a}} = \sum_{i=2b+c}^{2b+2c-2}d_{i+2b}x^i$,
and
\begin{equation*} 
B_{2_{a}} =  \sum_{i=b+c}^{2b+c-1}d_{i+2b}x^i +
  \sum_{i=b}^{2b+c-2}d_{i+2b+c}x^i +
  \sum_{i=c}^{b+2c-2}d_{i+3b}x^i +
  \sum_{i=0}^{b+c-2}d_{i+3b+c}x^i.
\end{equation*}

Proceeding with the reduction now of the second sum in
Equation (\ref{eq:A1}), we obtain
\begin{equation*} 
D_{2_{b}} =  \sum_{i=a}^{a+b-2}d_{i+2m-a-b}x^i +
  \sum_{i=b}^{2b-2}d_{i+2m-2b}x^i +
  \sum_{i=c}^{b+c-2}d_{i+2m-b-c}x^i +
  \sum_{i=0}^{b-2}d_{i+2m-b}x^i.
\end{equation*}
Clearly, $D_{2_b}$ is already reduced, and thus
$A_{2_{b}} = 0$,
and
\begin{equation*} 
B_{2_{b}} = \sum_{i=b+c}^{2b+c-2}d_{i+2b+c}x^i +
  \sum_{i=b}^{2b-2}d_{i+2b+2c}x^i +
  \sum_{i=c}^{b+c-2}d_{i+3b+c}x^i +
  \sum_{i=0}^{b-2}d_{i+3b+2c}x^i.
\end{equation*}

We finally reduce the third and last sum in
Equation~(\ref{eq:A1}):
\begin{equation*} 
D_{2_{c}} = \sum_{i=a}^{a+c-2}d_{i+2m-a-c}x^i +
  \sum_{i=b}^{b+c-2}d_{i+2m-b-c}x^i +
  \sum_{i=c}^{2c-2}d_{i+2m-2c}x^i +
  \sum_{i=0}^{c-2}d_{i+2m-c}x^i.
\end{equation*}
Again, we easily check that $D_{2_c}$ is reduced and so
$A_{2_{c}} = 0$,
and
\begin{equation*} 
B_{2_{c}} =  \sum_{i=b+c}^{b+2c-2}d_{i+3b}x^i +
  \sum_{i=b}^{b+c-2}d_{i+3b+c}x^i +
  \sum_{i=c}^{2c-2}d_{i+4b}x^i +
  \sum_{i=0}^{c-2}d_{i+4b+c}x^i.
\end{equation*}

\vspace{-0.1cm}

Concluding, $A_2$ is given by
\begin{equation} \label{eq:k3:A2}
  A_2 = A_{2_a} + A_{2_b} + A_{2_c}
      = \sum_{i=m}^{2a-2}d_{i+2m-2a}x^i,
\end{equation}
and
$B_2 = B_{2_a} + B_{2_b} + B_{2_c}$
is
\begin{equation} \label{eq:B2}
\begin{split}
  B_2 =
& \sum_{i=b+c}^{2b+c-1}d_{i+2b}x^i + \sum_{i=c}^{b+2c-2}d_{i+3b}x^i +  \sum_{i=b+c}^{b+2c-2}d_{i+3b}x^i + \sum_{i=c}^{2c-2}d_{i+4b}x^i +       \\
& \sum_{i=b}^{2b+c-2}d_{i+2b+c}x^i + \sum_{i=b+c}^{2b+c-2}d_{i+2b+c}x^i + \sum_{i=b}^{2b-2}d_{i+2b+2c}x^i  + \sum_{i=0}^{b+c-2}d_{i+3b+c}x^i +    \\
& \sum_{i=c}^{b+c-2}d_{i+3b+c}x^i  + \sum_{i=b}^{b+c-2}d_{i+3b+c}x^i +     \sum_{i=0}^{b-2}d_{i+3b+2c}x^i   + \sum_{i=0}^{c-2}d_{i+4b+c}x^i.
\end{split}
\end{equation}

\subsection{Determining $A_3$ and $B_3$}
Dividing the leading term of $A_2$ in Equation~(\ref{eq:k3:A2})
by $f$, we have
\begin{equation*}
D_{3} = \sum_{i=b+c}^{b+2c-2}d_{i+3b}x^{i} +
  \sum_{i=b}^{b+c-2}d_{i+3b+c}x^{i} +
  \sum_{i=c}^{2c-2}d_{i+4b}x^{i} +
  \sum_{i=0}^{c-2}d_{i+4b+c}x^{i}.
\end{equation*}
We have that $D_3$ is reduced so $A_3=0$ and
\begin{equation} \label{eq:B3}
B_{3} = \sum_{i=b+c}^{b+2c-2}d_{i+3b}x^{i} +
  \sum_{i=b}^{b+c-2}d_{i+3b+c}x^{i} +
  \sum_{i=c}^{2c-2}d_{i+4b}x^{i} +
  \sum_{i=0}^{c-2}d_{i+4b+c}x^{i}.
\end{equation}

\subsection{The number of terms in $A_r$ and $B_r$}

Let $G(i) = 1$ if $i > 0$ and $G(i) = 0$ if $i \le 0$. Let $r$ be a reduction step. It is clear now that the precise number of terms for $A_r$ and $B_r$, for $r \ge 0$, can be obtained using $k_{b+c}$, $k_b$ and $k_c$ given in Equations (\ref{eq:numero:reducoes}), (\ref{eq:numero:redb}) and (\ref{eq:numero:redc}). We have:
\begin{enumerate}[i)]
\item The number of terms of $A_0$ and $B_0$ is $1$.
\item For $r >0$, the number of terms of $A_r$ is $G(k_{b+c} - r) + G(k_b -r) + G(k_c -r)$, while the number of terms of $B_r$ is $4$ times the number of terms of $A_{r-1}$.
\end{enumerate}

\section{The family of polynomials $f(x) = x^{2b+1} +
x^{b+1} + x^b + x + 1$} \label{twored}

In this section, we consider the case when $c=1$, that is, when $k_{b+c}=2$,
as given in Equation~(\ref{eq:numero:reducoes}). The polynomials in this
subfamily have the form $f(x) = x^{2b+1} + x^{b+1} + x^b + x + 1$. For the
subfamily treated in this section, since $k_{b+c}=2$, we immediately get
$A_2 = 0$ and the expressions in the previous section simplify. As a
consequence, the desired reduction is given by
\begin{equation*} 
  D_{red} = B_0 + B_1 + B_2.
\end{equation*}

Using Equations~(\ref{eq:B0}),~(\ref{eq:B1})~and~(\ref{eq:B2}), we obtain
\begin{equation} \label{eq:B012:c1}
\begin{split}
  D_{red}
& =  \sum_{i=0}^{2b}d_ix^i + \sum_{i=b+1}^{2b}d_{i+b}x^i +  \sum_{i=1}^{b}d_{i+2b}x^i + \sum_{i=1}^{b}d_{i+3b}x^i +  \sum_{i=b}^{2b}d_{i+b+1}x^i + \\
& \sum_{i=0}^{b-1}d_{i+2b+1}x^i +  \sum_{i=b+1}^{2b-1}d_{i+2b+1}x^i + \sum_{i=b}^{2b-2}d_{i+2b+2}x^i + \sum_{i=0}^{b-2}d_{i+3b+2}x^i + d_{3b+1}.
\end{split}
\end{equation}

\begin{figure*}[hbt]
   \centering
   \includegraphics[width = .85\textwidth]{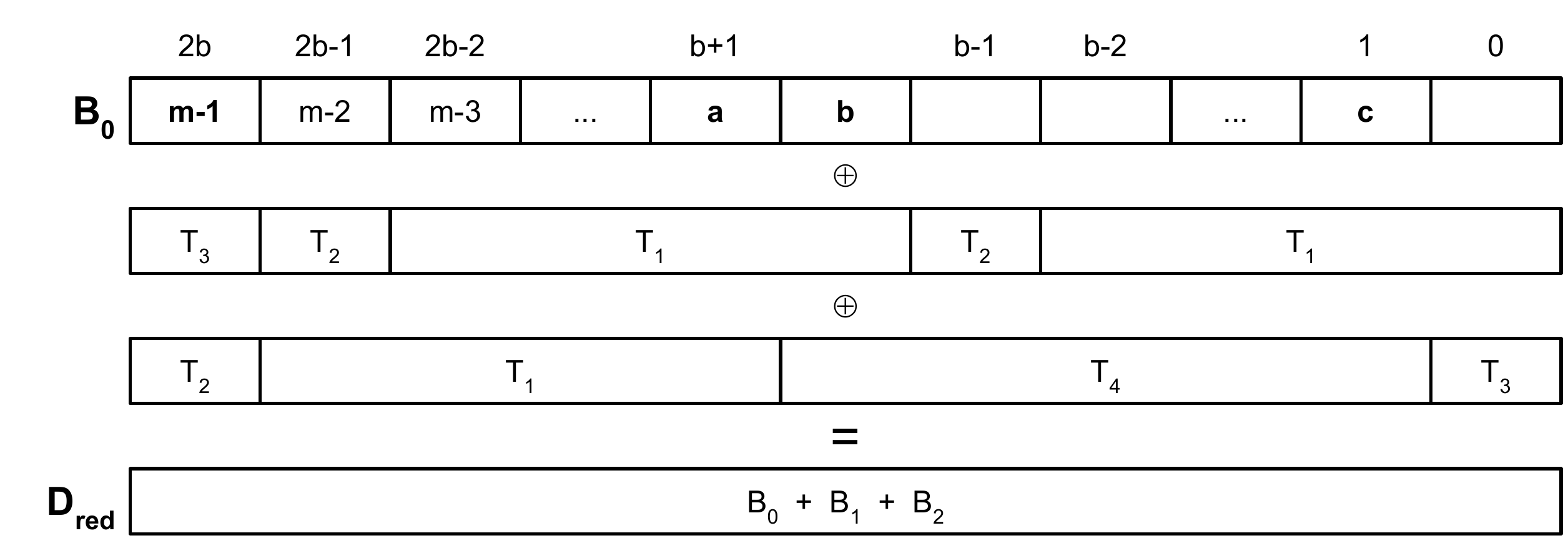}
   \caption{Representation of the reduction by $f(x)=x^{2b+1} + x^{b+1} + x^b + x + 1$.}
   \label{fig:Fig-DRed-ka2}
\end{figure*}

A crucial issue that allows us to give improved results for
our family of pentanomials is the fact that redundancies
occur for $D_{red}$ in Equation~(\ref{eq:B012:c1}).
Let
\begin{align}
  T_1(j) & = \sum_{i=0}^{b-2}(d_{i+2b+1} + d_{i+3b+2})x^{i+j}, \nonumber
             \qquad 
  T_2(j) = d_{3b}x^{j},  \nonumber \\
  T_3(j) & = d_{3b+1}x^{j}, \qquad  \nonumber
  T_4(j) = \sum_{i=0}^{b-1}(d_{i+2b+1} + d_{i+3b+1})x^{i+j}. \nonumber
\end{align}%
A careful analysis of Equation~(\ref{eq:B012:c1}) reveals that $T_1$, $T_2$
and $T_3$ are used more than once, and hence, savings can occur.
We rewrite Equation~(\ref{eq:B012:c1}) as
\begin{equation} \label{eq:Dred:c1}
  \begin{split}
    D_{red} = &
    B_0 + T_1(0) + T_1(b) + T_1(b+1) + T_2(b-1) +\\
    & T_2(2b-1) + T_2(2b) + T_3(0) + T_3(2b) + T_4(1).
  \end{split}
\end{equation}
One can check that by plugging $T_1$, $T_2$, $T_3$ and $T_4$ in Equation~(\ref{eq:Dred:c1}) we recover Equation~(\ref{eq:B012:c1}). Figure~\ref{fig:Fig-DRed-ka2} shows these operations. We remark that even though the first row in this figure is $B_0$, the following two rows are not $B_1$ and $B_2$. Indeed, those rows are obtained from $B_1$ and $B_2$ together with the optimizations provided by $T_1$, $T_2$, $T_3$ and $T_4$. 

Using Equation~(\ref{eq:Dred:c1}), the number $N_{\oplus}$ of XOR operations is
\[ N_{\oplus} = 6b + 1 = 3m -2. \]
It is also easy to see from Figure~\ref{fig:Fig-DRed-ka2} that the time delay is $3T_X$, where $T_X$ is the delay of one 2-input XOR gate.

We are now ready to provide Algorithm~\ref{alg:ka2} for computing $D_{red}$ given in Equation~(\ref{eq:Dred:c1}), and as explained in Figure~\ref{fig:Fig-DRed-ka2}, for the pentanomials $f(x) = x^{2b+1} + x^{b+1} + x^b + x + 1$.
\begin{algorithm}[htb] \small
\SetKwData{Left}{left}\SetKwData{This}{this}\SetKwData{Up}{up}
\SetKwFunction{Union}{Union}\SetKwFunction{FindCompress}{FindCompress}
\SetKwInOut{Input}{input}\SetKwInOut{Output}{output}
 \Input{$D_{0} = d[4b \ldots 0]$ bits vector of length $4b+1$}
 \Output{$D_{red}$}
 \For{$i \gets 0$ \KwTo $b-2$}{
    $T_1[i] \gets  d[i+2b+1] \oplus d[i+3b+2]$;\Comment{Definition of $T_1$}
 }
 \For{$i \gets 0$ \KwTo $b-1$}{
    $T_4[i] \gets  d[i+2b+1] \oplus d[i+3b+1]$;\Comment{Definition of $T_4$}
 }
 $D_{red}[0] \gets d[0] \oplus T_1[0] \oplus d[3b+1]$; \Comment{Column 0 of Fig.~\ref{fig:Fig-DRed-ka2}}

 \For{$i \gets 1$ \KwTo $b-2$}{
    $D_{red}[i] \gets  d[i] \oplus T_1[i] \oplus T_4[i-1]$; \Comment{Columns $1$ to $b-2$ of Fig.~\ref{fig:Fig-DRed-ka2}}
 }
 $D_{red}[b-1] \gets d[b-1] \oplus d[3b] \oplus T_4[b-2]$\;

 $D_{red}[b] \gets d[b] \oplus T_1[0] \oplus T_4[b-1]$\;

 \For{$i \gets b+1$ \KwTo $2b-2$}{
    $D_{red}[i] \gets  d[i] \oplus T_1[i-b] \oplus T_1[i-b-1]$; \Comment{Columns $b+1$ to $2b-2$ of Fig.~\ref{fig:Fig-DRed-ka2}}
 }
 $D_{red}[2b-1] \gets d[2b-1] \oplus d[3b] \oplus T_1[b-2]$\;

 $D_{red}[2b] \gets d[2b] \oplus d[3b+1] \oplus d[3b]$\;

  \Return $D_{red}$\;
 \caption{Computing $D_{red}$ when
$f(x) = x^{2b+1} + x^{b+1} + x^b + x + 1$.}
\label{alg:ka2}
\end{algorithm}

Putting all pieces together, we give next the main result of this section.

\begin{thm}
Algorithm~\ref{alg:ka2} correctly gives the reduction of a polynomial of degree at most $2m-2$ over $\ftwo$ by $f(x) = x^{2b+1} + x^{b+1} + x^b + x + 1$ involving $N_{\oplus} = 3m -2 = 6b + 1$ number of XORs operations and a time delay of $3T_X$.
\end{thm}

\section{Family $f(x) = x^{2b+c} + x^{b+c} + x^b + x^c + 1, c>1$}
\label{threered}
For polynomials of the form $f(x) = x^{2b+c} + x^{b+c} + x^b + x^c + 1$, $c>1$, we have that $k_{b+c}=3$, implying that $A_3 = 0$. The reduction is given by
\begin{equation*} 
  D_{red} = B_0 + B_1 + B_2 + B_3.
\end{equation*}
Using Equations~(\ref{eq:B0}),~(\ref{eq:B1}),~(\ref{eq:B2})
and ~(\ref{eq:B3}), we have that $D_{red}$ satisfies
\begin{equation}\label{eq:otm:ka3}
\begin{split}
D_{red}  =
& \sum_{i=0}^{2b+c-1}d_ix^i +
  \sum_{i=b+c}^{2b+c-1}d_{i+b}x^i + \sum_{i=c}^{b+c-1}d_{i+2b}x^i + \sum_{i=c}^{b+2c-2}d_{i+3b}x^i +   \\
&  \sum_{i=b}^{2b+c-1}d_{i+b+c}x^i +
 \sum_{i=0}^{b-1}d_{i+2b+c}x^i + \sum_{i=b+c}^{2b+c-2}d_{i+2b+c}x^i + \\
& \sum_{i=b}^{2b-2}d_{i+2b+2c}x^i + \sum_{i=0}^{c-1}d_{i+3b+c}x^i +
  \sum_{i=0}^{b-2}d_{i+3b+2c}x^i.
\end{split}
\end{equation}

\begin{figure*}[htb]
   \centering
   \includegraphics[width = \textwidth]{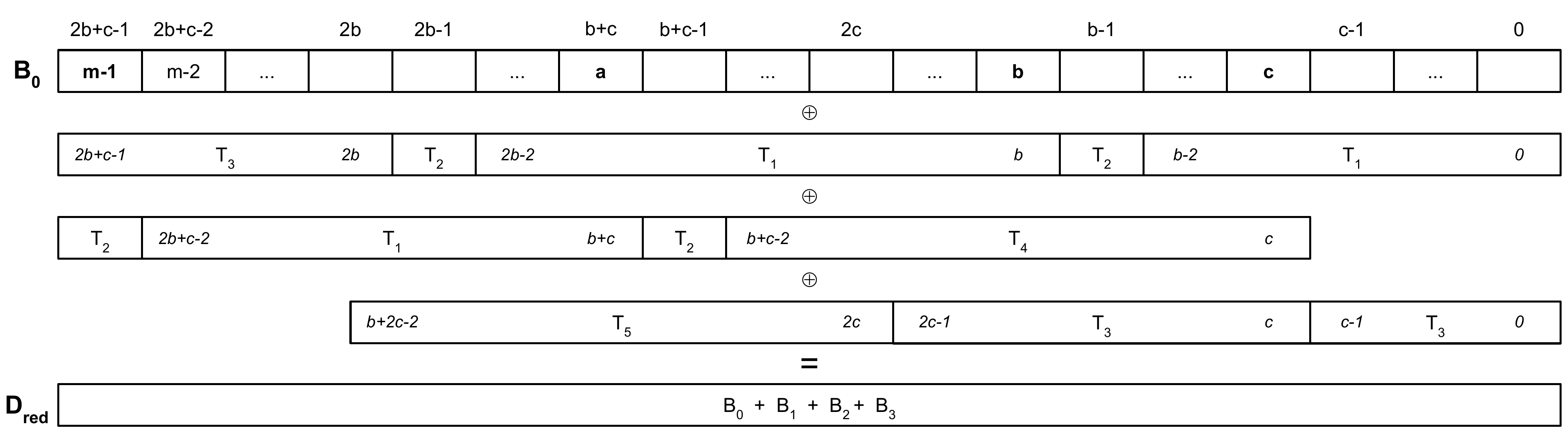}
   \caption{Representation of the reduction by $f(x) = x^{2b+c} + x^{b+c} + x^b + x^c + 1, c>1$.}
   \label{fig:Fig-DRed-ka3}
\end{figure*}

Let
\begin{align}
  T_{1}(j) & = \sum_{i=0}^{b-2}(d_{i+2b+c} + d_{i+3b+2c})x^{i+j}, \qquad \nonumber
  T_2(j) =  d_{3b+c-1}x^j,                                    \nonumber \\
  T_3(j)   & = \sum_{i=0}^{c-1}d_{i+3b+c}x^{i+j}, \quad \nonumber T_4(j) = \sum_{i=0}^{b-2}d_{i+2b+c}x^{i+j},                 \quad \nonumber T_5(j) = \sum_{i=0}^{b-2}d_{i+3b+2c}x^{i+j}.                \nonumber
\end{align}
Again, a careful analysis of
Equation~(\ref{eq:otm:ka3}) shows that $T_1$, $T_2$ and $T_3$ are
used more than once. Thus, we can rewrite Equation~(\ref{eq:otm:ka3})
for $D_{red}$ as
\begin{equation} \label{eq:Dred:clarger1}
\begin{split}
D_{red} = & B_0 + T_1(0) + T_1(b) + T_1(b+c) + \\
          & T_2(b-1) + T_2(b+c-1) + T_2(2b-1) + T_2(2b+c-1) + \\
          & T_3(0) + T_3(c) + T_3(2b) + T_4(c) +T_5(2c).
\end{split}
\end{equation}
Figure~\ref{fig:Fig-DRed-ka3} depicts these operations.
Using Equation~(\ref{eq:Dred:clarger1}) and 
Figure~\ref{fig:Fig-DRed-ka3}, we have Algorithm~\ref{alg:ka3}. For code efficiency reasons, in contrast to Algorithm~\ref{alg:ka2}, in Algorithm~\ref{alg:ka3} we separate the last line before the equality in Figure~\ref{fig:Fig-DRed-ka3}. The additions of this last line are done in lines $17$ to $20$ of Algorithm~\ref{alg:ka3}. As a consequence, lines $3$ to $16$ of Algorithm~\ref{alg:ka3} include only the additions per column from $0$ to $2b+c-1$ of the first three lines in Figure~\ref{fig:Fig-DRed-ka3}.

\begin{algorithm}[ht] \small
\SetKwData{Left}{left}\SetKwData{This}{this}\SetKwData{Up}{up}
\SetKwFunction{Union}{Union}\SetKwFunction{FindCompress}{FindCompress}
\SetKwInOut{Input}{input}\SetKwInOut{Output}{output}
 \Input{$D_{0} = d[2b+c-1 \ldots 0]$ bits vector of length $2b+c$}
 \Output{$D_{red}$}
 \For{$i \gets 0$ \KwTo $b-2$}{
    $T_1[i] \gets  d[i+2b+1] \oplus d[i+3b+2c]$; \Comment{Definition of $T_1$}
 }
 \For{$i \gets 0$ \KwTo $c-1$}{
    $D_{red}[i] \gets  d[i] \oplus T_1[i]$; \Comment{Columns $0$ to $c-1$ of the first three lines of  Fig.~\ref{fig:Fig-DRed-ka3}}
 }
 \For{$i \gets c$ \KwTo $b-2$}{
    $D_{red}[i] \gets  d[i] \oplus T_1[i] \oplus d[i+2b]$\;
 }
 $D_{red}[b-1] \gets  d[b-1] \oplus  d[3b+c-1] \oplus d[3b-1]$\;
 \For{$i \gets b$ \KwTo $b+c-2$}{
    $D_{red}[i] \gets  d[i] \oplus T_1[i-b] \oplus d[i+2b]$\;
 }
 $D_{red}[b+c-1] \gets  d[b+c-1] \oplus  d[3b+c-1] \oplus T_1[c-1]$\;
 \For{$i \gets b+c$ \KwTo $2b-2$}{
    $D_{red}[i] \gets  d[i] \oplus T_1[i-b] \oplus T_1[i-b-c]$\;
 }
 $D_{red}[2b-1] \gets  d[2b-1] \oplus  d[3b+c-1] \oplus T_1[b-c-1]$\;
 \For{$i \gets 2b$ \KwTo $2b+c-2$}{
    $D_{red}[i] \gets  d[i] \oplus T_1[i-b-c] \oplus d[i+b+c]$\;
 }
 $D_{red}[2b+c-1] \gets  d[2b+c-1] \oplus  d[3b+c-1] \oplus d[3b-1]$\;
 \For{$i \gets 0$ \KwTo $c-1$}{
    $D_{red}[i] \gets  D_{red}[i] \oplus d[i+3b+c]$; \Comment{Columns $0$ to $c-1$ of the $4^{th}$ line of  Fig.~\ref{fig:Fig-DRed-ka3}}
 }
 \For{$i \gets c$ \KwTo $b+2c-2$}{
    $D_{red}[i] \gets  D_{red}[i] \oplus d[i+3b]$; \Comment{Cols $c$ to $b+2c-2$ of the $4^{th}$ line of  Fig.~\ref{fig:Fig-DRed-ka3}}
 }

 \Return $D_{red}$\;
 \caption{Computing $D_{red}$ for
$f(x) = x^{2b+c} + x^{b+c} + x^b + x^c + 1$.}
 \label{alg:ka3}
\end{algorithm}
\vspace{-.5cm}
The time delay is $3T_X$;
after removal of redundancies and not counting repeated
terms, we obtain that the number $N_\oplus$ of XORs is
\[ N_\oplus = 6b + 3c - 2 = 3m-2.\]

\begin{thm}
Algorithm~\ref{alg:ka3} correctly gives the reduction
of a polynomial of degree at most $2m-2$ over $\ftwo$
by $f(x) = x^{2b+c} + x^{b+c} + x^b + x^c + 1$ involving
$N_{\oplus} = 3m -2 = 6b + 3c-2$ number of XORs operations
and a time delay of $3T_X$.
\end{thm}

\subsection{Almost equally spaced pentanomials: the special case $b=2c$}
Consider the special case $b=2c$. In this case we obtain the almost equally spaced polynomials $f(x) = x^{5c} + x^{3c} + x^{2c} + x^c + 1$.
Our previous analysis when applied to these polynomials entails
\begin{equation}\label{eq:otm:ka3:b2c}
\begin{split}
D_{red} =
  & \sum_{i=0}^{5c-1}d_ix^i       + \sum_{i=3c}^{5c-1}d_{i+2c}x^i + \sum_{i=c}^{3c-1}d_{i+4c}x^i  + \sum_{i=c}^{4c-2}d_{i+6c}x^i  + \sum_{i=2c}^{5c-1}d_{i+3c}x^i + \\
  & \sum_{i=0}^{2c-1}d_{i+5c}x^i  + \sum_{i=3c}^{5c-2}d_{i+5c}x^i + \sum_{i=2c}^{4c-2}d_{i+6c}x^i + \sum_{i=0}^{c-1}d_{i+7c}x^i   + \sum_{i=0}^{2c-2}d_{i+8c}x^i.
\end{split}
\end{equation}

Let
\begin{align}
  T_{1}(j) & = \sum_{i=c}^{2c-2}(d_{i+5c} + d_{i+4c})x^{i+j}, \quad \nonumber
  T_2(j)    = \sum_{i=c}^{2c-2}(d_{i+8c} + d_{i+6c})x^{i+j}, \quad \nonumber \\
  T_3(j)   & =  d_{8c-1}x^j,                                  \quad \nonumber
  T_4(j)    = \sum_{i=0}^{c-1}d_{i+8c}x^{i+j},               \quad \nonumber
  T_5(j)    = \sum_{i=0}^{c-1}d_{i+5c}x^{i+j},               \quad \nonumber \\
  T_6(j)   & = \sum_{i=0}^{c-2}d_{i+7c}x^{i+j},               \quad \nonumber
  T_7(j)    = \sum_{i=4c}^{5c-1}d_{i+2c}x^{i+j}.                   \nonumber
\end{align}

\begin{figure*}[htb]
   \centering
   \includegraphics[width = \textwidth]{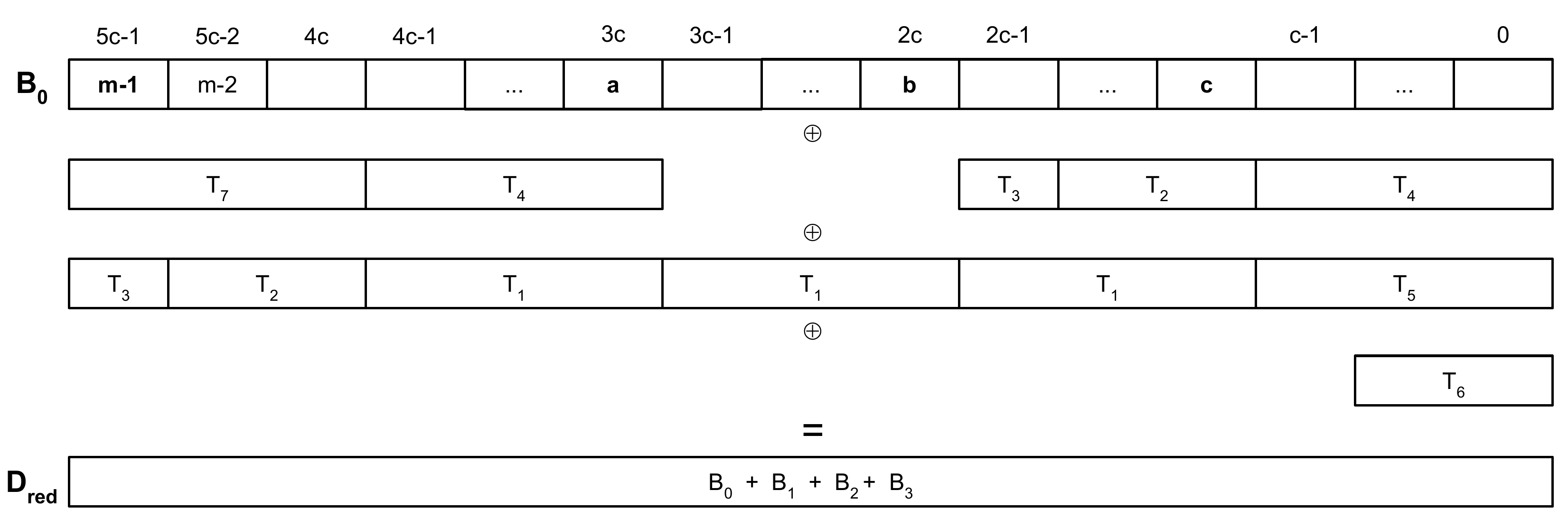}
   \caption{Representation of the reduction by the almost equally spaced pentanomials (the special case $b=2c$).}
   \label{fig:Fig-DRed-ka3-b-2c}
\end{figure*}

In the computation of $D_{red}$, $T_1$, $T_2$, $T_3$ and $T_4$ are used more than once.
Figure~\ref{fig:Fig-DRed-ka3-b-2c} shows, graphically, these operations. After removal of redundancies, the number $N_\oplus$ of XORs is
$ N_\oplus = 12c - 1 = \displaystyle\frac{12}{5}m -1.$
This number of XORs is close to $2.4m$ providing a saving of about $20\%$ with respect to the other pentanomials in our family. Irreducible pentanomials of this form are rare but they do exist,
for example, for degrees $5$, $155$ and $4805$. We observe that the extension $155$ is used in \cite{Vanstone}.

Using Equation~(\ref{eq:otm:ka3:b2c}) and
Figure~\ref{fig:Fig-DRed-ka3-b-2c}, we naturally have Algorithm~\ref{alg:b2c}.

\begin{algorithm}[htb] \small
\SetKwData{Left}{left}\SetKwData{This}{this}\SetKwData{Up}{up}
\SetKwFunction{Union}{Union}\SetKwFunction{FindCompress}{FindCompress}
\SetKwInOut{Input}{input}\SetKwInOut{Output}{output}
 \Input{$D_{0} = d[5c-1 \ldots 0]$ bits vector of length $5c$}
 \Output{$D_{red}$}
 \For{$i \gets 0$ \KwTo $c-2$}{
    $T_1[i] \gets  d[i+6c] \oplus d[i+5c]$\; \Comment{Definition of $T_1$}

 }
 \For{$i \gets 0$ \KwTo $c-2$}{
    $T_2[i] \gets  d[i+9c] \oplus d[i+7c]$\;\Comment{Definition of $T_2$}

 }
 \For{$i \gets 0$ \KwTo $c-2$}{
    $D_{red}[i] \gets  d[i] \oplus d[i+8c] \oplus d[i+5c] \oplus d[i+7c]$\;
 }
 $D_{red}[c-1] \gets  d[c-1] \oplus d[9c-1] \oplus d[6c-1]$\;

 \For{$i \gets c$ \KwTo $2c-2$}{
    $D_{red}[i] \gets  d[i] \oplus T_1[i-c] \oplus T_2[i-c]$\;
 }
 $D_{red}[2c-1] \gets d[2c-1]  \oplus d[8c-1] \oplus T_1[c-1]$\;

 \For{$i \gets 2c$ \KwTo $3c-1$}{
    $D_{red}[i] \gets  d[i] \oplus T_1[i-2c]$\;
 }
  \For{$i \gets 3c$ \KwTo $4c-1$}{
    $D_{red}[i] \gets  d[i] \oplus T_1[i-3c] \oplus d[i+5c]$\;
 }
 \For{$i \gets 4c$ \KwTo $5c-2$}{
    $D_{red}[i] \gets  d[i] \oplus T_2[i-4c] \oplus d[i+2c]$\;
 }
 $D_{red}[5c-1] \gets d[5c-1]  \oplus d[8c-1]  \oplus d[7c-1]$\;

  \Return $D_{red}$\;
 \caption{Computing $D_{red}$ for
$f(x) = x^{5c} + x^{3c} + x^{2c} + x^c + 1$.}
 \label{alg:b2c}
\end{algorithm}

\section{Multiplier in $\mathbb{F}_{2}[x]$, complexity analysis and comparison}\label{sec:comparisons}

So far, we have focused on the second step of the algorithm,
that is, on the reduction part. For the first step, the
multiplication part, we simply use a standard Karatsuba recursive algorithm
implementation; see Algorithm~\ref{koa:poly}.

Recursivity in hardware can be an issue; see \cite{von2005efficient} and \cite{machhout2008efficient},
for example, for efficient hardware implementations of polynomial multiplication in finite fields
using Karatsuba's type algorithms.
\begin{algorithm} \small
	\SetAlgoLined
	\DontPrintSemicolon
  \SetKwData{Left}{left}\SetKwData{This}{this}\SetKwData{Up}{up}
  \SetKwFunction{Union}{Union}\SetKwFunction{FindCompress}{FindCompress}
  \SetKwInOut{Input}{input}\SetKwInOut{Output}{output}
	\Input{$A(x) = \sum_{i = 0}^{m-1} a_ix^i$ and  $B(x) = \sum_{i = 0}^{m-1} b_ix^i$}
	\Output{$C(x) = A(x)B(x) = \sum_{i = 0}^{2m-2} c_ix^i$ }

	\SetKwFunction{FMain}{Karatsuba}
	\SetKwProg{Fn}{Function}{:}{}
	\Fn{\FMain{$A$, $B$}}{

		$m \gets maxDegree(A, B)$\; \Comment{compute the larger degree between polynomials $A$ and $B$ }

		\If{m = 0}{

			\Return $(A$ \& $B)$\;\Comment{\& is a bitwise AND operator}
		}
		$m2 = floor(m/2)$\; \Comment{split $A$ and $B$}

		$\text{high}_{a}, \text{low}_{a} \gets split(A,m2)$\;

		$\text{high}_{b}, \text{low}_{b} \gets split(B,m2)$\;

		$d_0 \gets \text{Karatsuba}(\text{low}_a, \text{low}_b)$\; \Comment{recursive call of Karatsuba}

		$d_1 \gets \text{Karatsuba}((\text{low}_a\oplus\text{high}_a), (\text{low}_b \oplus \text{high}_b))$\; \Comment{recursive call of Karatsuba}

		$d_2 \gets \text{Karatsuba}(\text{high}_a, \text{high}_b)$\; \Comment{recursive call of Karatsuba}

		$c \gets d_2x^{m} \oplus (d_1 \oplus d_2 \oplus d_0)x^{m2} \oplus d_0$\;

		\Return $c$\;
	}
	\textbf{End Function}
	\caption{Karatsuba Algorithm for $\ftwom$}
	\label{koa:poly}
\end{algorithm}

As can be seen our multiplier consists of two steps. The first is the multiplication itself using Karatsuba arithmetic or, if necessary, the school book method, and the second is the reduction described in the 
previous sections. The choice of the first step method will
basically depend on whether the application requirement is 
to minimize area (Karatsuba), i.e., the number of ANDs and XORs gates,
or to minimize the arithmetic delay (School book); see~\cite{FanHasan} for 
several variants of both the schoolbook and Karatsuba 
algorithms. Minimizing the area is interesting in applications 
that need to save power at the expense of a longer runtime.

We chose the Karatsuba multiplier since our goal is to minimize 
the area, i.e. to minimize the number of gates AND and XOR. 
A summary of our results compared with related works is given 
in Tables~\ref{table:costs:pentanomials} and~\ref{tab:1}. 
Table~\ref{table:costs:pentanomials} presents comparison 
costs among multipliers that perform two steps for the 
multiplication, that is, they execute a multiplication 
followed by a reduction. The table shows the multiplication 
algorithm used in each case. Table~\ref{tab:1} gives a 
comparison among the state-of-the-art bit multipliers in the 
literature. The main target for us is \cite{park2014new} 
since it presents the smallest area in the literature. 
However, Type $3$ polynomials are also considered; this is 
another practically relevant family of polynomials. With 
respect to Karatsuba variants, Table~3 of survey~\cite{FanHasan} 
shows asymptotic complexities of several Karatsuba 
multiplication algorithms without reduction.

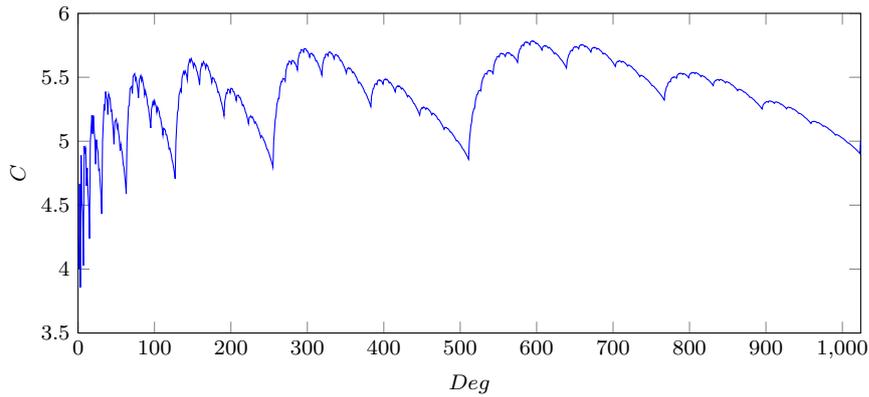
\begin{figure}
\centering
\begin{tikzpicture}
    \begin{axis}[
      ticklabel style = {font=\small}, xlabel={$Deg$},
      ylabel={$C$},ymin=3.5,ymax=6,xmin=0,xmax=1024,width=\textwidth, height=.3\textheight,ytick={3.5,4,4.5,5,5.5,6}]
      \addplot[color=blue] table[x=Deg,y=C] {c-karatsuba.txt};
    \end{axis}
\end{tikzpicture}
\caption{Karatsuba constant for degrees up to $1024$.}
\label{fig:ka:cte}
\end{figure}

For each entry in Table~\ref{table:costs:pentanomials}, we 
give the multiplication algorithm and the amount 
of gates AND, XOR as well its delay. We point that for 
\cite{Masoleh:2004} and~\cite{wu1999}, their multipliers 
are general for any pentanomial with $a \leq \frac{m}{2}$ 
instead of for a specific family such as~\cite{RodriguezKoc}.
In the case of our family, in addition to the number of XORs 
for the reduction, we include the cost for the multiplication 
due to the recursive Karatsuba implementation multiplier, 
that is, the XOR count is formed by the sum of the XORs of 
the Karatsuba multiplier and the ones of the reduction part. 
In our implementation, the constant of Karatsuba is strictly 
less than $6$; see 
Figure~\ref{fig:ka:cte}
for degrees up to $1024$. As can be seen, for degrees powers 
of $2$ minus $1$ ($2^k -1, k \geq 1$), the constant achieves 
local minimum. For the number of AND gates, we provide an 
interval. The actual number of AND gates depends on the value 
of $m$; it only reaches a maximum when $m = 2^{k}-1$, for 
$k \geq 1$.

\begin{table}[htb]
\centering
\caption{Two steps multipliers cost comparison for different family of pentanomials.}
\label{table:costs:pentanomials}
{\small
\begin{adjustbox}{max width=\textwidth}
\begin{tabular}{l|c|c|c}
\hline 
\multicolumn{4}{c}{$x^m + x^a +x^b+x^c+1$~\cite{wu1999,RodriguezKoc}, Multiplication algorithm: Schoolbook.}\\ \hline
\multicolumn{1}{c|}{Costs} & \multicolumn{1}{c|}{\#AND} & \multicolumn{1}{c|}{\#XOR} & \multicolumn{1}{c}{Delay} \\ \hline
Reduction & 0 & $4(m-1)$ & $3T_X$ \\
Multiplication & $m^2$ & $(m-1)^2$ & $T_A + (\lceil \log_2{m} \rceil ) T_X$ \\
Multiplier & $m^2$ & $m^2 + 2m -3 $ & $T_A + (3 + \lceil \log_2{m} \rceil ) T_X$ \\ \hline 

\multicolumn{4}{c}{Type $I$ - $x^m + x^{n+1} +x^n+x+1$~\cite{RodriguezKoc}, Multiplication algorithm: Mastrovito-like Multiplier.}\\ \hline
\multicolumn{1}{c|}{Costs} & \multicolumn{1}{c|}{\#AND} & \multicolumn{1}{c|}{\#XOR} & \multicolumn{1}{c}{Delay} \\ \hline
Reduction & 0 & $3m+2n-1$ & $3T_X$\\
Multiplication & $m^2$ & $ m^2 - 2m + 1 $&  $T_A + (\lceil \log_2{m} \rceil ) T_X$\\
Multiplier & $m^2$ & $m^2 + m + 2n$ & $T_A + (3 + \lceil \log_2{m} \rceil ) T_X$ \\ \hline

\multicolumn{4}{c}{Type $I$ - $x^m + x^{n+1} +x^n+x+1$~\cite{Masoleh:2004}, Multiplication algorithm: Mastrovito-like Multiplier.}\\ \hline
\multicolumn{1}{c|}{Costs} & \multicolumn{1}{c|}{\#AND} & \multicolumn{1}{c|}{\#XOR} & \multicolumn{1}{c}{Delay} \\ \hline
Reduction &$0$ & $3m-2$ & $3T_X$ \\
Multiplication & $m^2$ & $m^2-2m+1$ & $T_A + (\lceil \log_2{(m-1)} \rceil)T_X $ \\
Multiplier & $m^2$ & $m^2+m ^\dagger$ & $T_A + (3+\lceil \log_2{(m-1)} \rceil)T_X  $ \\ \hline

\multicolumn{4}{c}{Type $II$ - $x^m + x^{n+2} +x^{n+1}+x^n+1$~\cite{RodriguezKoc}, Multiplication algorithm: Dual basis.}\\ \hline
\multicolumn{1}{c|}{Costs} & \multicolumn{1}{c|}{\#AND} & \multicolumn{1}{c|}{\#XOR} & \multicolumn{1}{c}{Delay} \\ \hline
Reduction   & $0$ & $3m- \lceil (m-2)/2 \rceil +3n - 4$ & $3T_X$ \\
Multiplication  & $m^2$ & $m^2-m$ & $T_A + (\lceil \log_2{m} \rceil ) T_X$ \\
Multiplier & $m^2$ & $m^2 + 2m- \lceil (m-2)/2 \rceil +3n - 4$ & $T_A + (3 + \lceil \log_2{m} \rceil ) T_X$ \\ \hline

\multicolumn{4}{c}{$x^m + x^a +x^b+x^c+1, c > 1$~\cite{Masoleh:2004}, Multiplication algorithm: Mastrovito-like Multiplier.}\\ \hline
\multicolumn{1}{c|}{Costs} & \multicolumn{1}{c|}{\#AND} & \multicolumn{1}{c|}{\#XOR} & \multicolumn{1}{c}{Delay} \\ \hline
Reduction &$0$ & $4m -4$ & $4T_X$ \\
Multiplication & $m^2$ & $m^2-2m+1$ & $T_A + (\lceil \log_2{(m-1)} \rceil)T_X$ \\
Multiplier & $m^2$ & $m^2+2m-3$ & $T_A + (4+\lceil \log_2{(m-1)} \rceil)T_X $ \\ \hline 

\multicolumn{4}{c}{Ours - $x^{2b+c} + x^{b+c} +x^b+x^c+1$, Multiplication algorithm: Karatsuba.}\\ \hline
\multicolumn{1}{c|}{Costs} & \multicolumn{1}{c|}{\#AND} & \multicolumn{1}{c|}{\#XOR} & \multicolumn{1}{c}{Delay} \\ \hline
Reduction & 0 & $3m-2$ & $3T_X$\\
Multiplication & $( 3^{\lfloor \log_2{m} \rfloor}, 3^{\lfloor \log_2{m} \rfloor + 1}]$ & $< 6m^{\log_2{3}}$ & $T_A + 3\lceil \log_2{(m-1)} \rceil T_X$\\
Multiplier & $( 3^{\lfloor \log_2{m} \rfloor}, 3^{\lfloor \log_2{m} \rfloor + 1}]$ & $< 6m^{\log_2{3}} + 3m-2$ & $T_A + 3 (\lceil \log_2{(m-1)} \rceil+1)T_X$ \\ \hline 

\multicolumn{4}{c}{Ours - $x^{5c} + x^{3c} +x^{2c}+x^c+1$, Multiplication algorithm: Karatsuba.}\\ \hline
\multicolumn{1}{c|}{Costs} & \multicolumn{1}{c|}{\#AND} & \multicolumn{1}{c|}{\#XOR} & \multicolumn{1}{c}{Delay} \\ \hline
Reduction & 0 & $(12/5)m - 1$ & $3T_X$\\
Multiplication & $( 3^{\lfloor \log_2{m} \rfloor}, 3^{\lfloor \log_2{m} \rfloor + 1}]$ & $< 6m^{\log_2{3}}$ & $T_A + 3\lceil \log_2{(m-1)} \rceil T_X$ \\
Multiplier & $( 3^{\lfloor \log_2{m} \rfloor}, 3^{\lfloor \log_2{m} \rfloor + 1}]$ & $< 6m^{\log_2{3}} + (12/5)m - 1$ & $T_A + 3 (\lceil \log_2{(m-1)} \rceil + 1)T_X$\\ \hline 

\end{tabular}
\end{adjustbox}
{\scriptsize $\dagger$ There
is an additional XOR to reduce the time delay; see~\cite[page 955]{Masoleh:2004}.}
}
\end{table}

In Table~\ref{tab:1}, we provide the number of XORs and 
ANDs gates for Type $1$ and Type $2$ families in
\cite{park2014new} and \cite{RodriguezKoc}, 
Type $3$ in \cite{Masoleh:2004} and our family of 
pentanomials. We point out that in~\cite{park2014new} 
the authors compute multiplication and reduction as a 
unique block with a divide-and-conquer approach using
squaring. In contrast, we separate these two parts and 
use Karatsuba for the multiplier followed by our 
reduction algorithm. 

\begin{table}[htb]
\centering
\caption{Space and time complexities of state-of-the-art bit 
multipliers.}
\label{tab:1}
 {\small
\begin{adjustbox}{max width=\textwidth}
\begin{tabular}{lllll}
Type & & \multicolumn{1}{c}{\# XOR} & \# AND & Delay \\ \hline

Type 1 & \multicolumn{3}{l}{$x^m + x^{b+1} + x^{b} + x + 1$,  $1 < b \leq \frac{m}{2}-1$} &  \\

\small{\cite{park2014new}~~} & $b$ is odd & $\displaystyle\frac{3m^2+24m+8b+21}{4}$ & $\displaystyle\frac{3m^2+2m-1}{4}$ & $T_A + (3+\lceil \log_2 (m+1)\rceil)T_x$ \\

\small{\cite{park2014new}~~} & $b$ is even & $\displaystyle\frac{3m^2+24m+8b+17}{4}$ & $\displaystyle\frac{3m^2+2m-1}{4}$  & $T_A + (3+\lceil \log_2 (m+1)\rceil)T_x$\\ \hline

Type 2 & \multicolumn{3}{l}{$x^m + x^{c+2} + x^{c+1} + x^{c} + 1$} & \\

\small{\cite{park2014new}~~} & $c$ is odd, $c \leq \frac{3}{8}(m-7)$ & $\displaystyle\frac{3m^2+24m+14c+35}{4}$ & $\displaystyle\frac{3m^2+2m-1}{4}$ & $T_A + (3+\lceil \log (m+1)\rceil)T_x$ \\

\small{\cite{park2014new}~~} & $c$ is even, $c \leq \frac{m}{2}-1$ & $\displaystyle\frac{3m^2+24m+14c+45}{4}$ & $\displaystyle\frac{3m^2+2m-1}{4}$ & $T_A + (3+\lceil \log (m+1)\rceil)T_x$ \\

\small{\cite{RodriguezKoc}~~} &$c > 1$  & $\displaystyle{m^2+2m - \lceil (m-2)/2 \rceil +3n -4}$ & $\displaystyle{m^2}$ & $T_A + (3+\lceil \log (m-1)\rceil)T_x$ \\

\small{\cite{RodriguezKoc}~~} &$c = 1$  & $\displaystyle{m^2+m-2}$ & $\displaystyle{m^2}$ & $T_A + (3+\lceil \log_2 (m-1)\rceil)T_x$ \\ \hline

Type 3 & \multicolumn{3}{l}{$x^{m} + x^{m-c} + x^{m-2c} + x^{m-3c} + 1$ } & \\

\small{\cite{Masoleh:2004}~~} &$\frac{m-1}{4} \leq c \leq \frac{m-1}{3}$  & $m^2 + m - c -1$ & $\displaystyle{m^2}$ & $T_A + (3+\lceil \log_2 (m-1)\rceil)T_x$\\

\small{\cite{Masoleh:2004}~~} &$\frac{m-1}{5} \leq c < \frac{m-1}{4}$  & $m^2 + 2m - 5c - 2$ & $\displaystyle{m^2}$ &  $T_A + (3+\lceil \log_2 (m-1)\rceil)T_x$\\

\small{\cite{Masoleh:2004}~~} &$\frac{m-1}{8} \leq c < \frac{m-1}{5}$  & $m^2 + m - 2$ & $\displaystyle{m^2}$ & $T_A + (3+\lceil \log_2 (m-1)\rceil + 1)T_x$ \\ \hline

Ours & \multicolumn{3}{l}{$x^{2b+c}+x^{b+c}+x^b + x^c + 1$} & \\

\small{Ours} & $c\ge 1$, $b\neq 2c$ &  $< 6m^{\log_2{3}} + 3m-2$ & $( 3^{\lfloor \log_2{m} \rfloor}, 3^{\lfloor \log_2{m} \rfloor + 1}]$ & $T_A + 3 (\lceil \log_2{(m-1)} \rceil)T_x$ \\ 

\small{Ours} & $c\ge 1$, $b = 2c$ & $< 6m^{\log_2{3}} + \frac{12}{5}m - 1$ & $( 3^{\lfloor \log_2{m} \rfloor}, 3^{\lfloor \log_2{m} \rfloor + 1}]$ & $T_A + 3 (\lceil \log_2{(m-1)} \rceil + 1)T_x$\\ \hline

\end{tabular}
\end{adjustbox}
}
\end{table}

The costs for using our pentanomials for degrees proposed by NIST can be found in Table~\ref{table:fixed:results}. The amount of XOR and AND gates are the exact value obtained from Table~\ref{table:costs:pentanomials}. The delay costs can be separated in $T_A$ and $T_X$, delay for AND gates and XOR gates, respectively. The delay for AND gates is due to only $1$ AND gate at the lowest level of the Karatsuba recursion.
The delay for the XOR gates in the Karatsuba multiplier is $3 \lceil \log_2{(m-1)} \rceil$ since there are $3$ delay XORs per
level of the Karatsuba recursion. For the reduction part, we only have $3$ delay XORs. Hence, the total number of XOR delays is $3 \lceil \log_2{(m-1)} \rceil + 3$.

Table~\ref{tab:2} shows the number of irreducible pentanomials of
degrees $163$, $283$ and $571$ for the families considered since those
are NIST degrees where pentanomials have been recommended~\cite{NIST}.
Analyzing the table, we have that family Type $1$ has the most irreducible
pentanomials, but few of them have degrees recommended by NIST \cite{NIST}. The first family of Type 2, proposed in \cite{park2014new}, has restrictions in the range of $c$; this family presents the highest number of representatives with NIST degrees of interest. The second family of Type 2, proposed in \cite{RodriguezKoc}, has no restrictions for $c$; this family presents the largest number of irreducible polynomials. Type $3$ is the special case from \cite{Masoleh:2004}. Our family for $b \not = 2c$ has less irreducible polynomials and it has no irreducible polynomials with degrees $163$, $283$ and $571$. In the other side, when $b \neq 2c$ our family has $730$ polynomials of degrees up to $1024$ and it presents $5$ pentanomials of NIST degrees.

In the following we comment on the density of irreducible pentanomials in our family.
Table~\ref{tab:penta:irredut:2} lists all irreducible pentanomials of our family for degrees up
to $1024$; $N_\oplus$ is the number of XORs required for the reduction.
We leave as an open problem to mathematically characterize under which
conditions our pentanomials are irreducible.

\begin{table}[htb]
\centering
\caption{Costs for fixed degree pentanomials proposed by NIST.}
\label{table:fixed:results}
{\small
\begin{tabular}{c|r|r|r|r|r}
\hline
\multirow{2}{*}{Degree} & \multicolumn{3}{c|}{XORs}     & \multirow{2}{*}{ANDs} & \multicolumn{1}{c}{\multirow{2}{*}{Delay}} \\ \cline{2-4}

& Karatsuba & Reduction & Total &  & \multicolumn{1}{c}{} \\ \hline

$163$ & $17,944$ & $487$ & $18,431$ & $4,419$ & $T_A + 27T_X$ \\ \hline

$283$ & $43,162$ & $847$ & $44,009$ & $10,305$ & $T_A + 30T_X$ \\ \hline

$571$ & $132,280$ & $1,711$ & $133,991$ & $31,203$ &$T_A + 33T_X$ \\ \hline
\end{tabular}
}
\end{table}

\begin{table}[!htb]
\centering
\caption{Number of irreducible pentanomials for NIST degrees.}
\label{tab:2}
\begin{tabular}{l|r|r|r|r}
  Type & \#Irred. & $163$ & $283$ & $571$   \\ \hline
Type $1$~\cite{park2014new}  & $2025$ & $1$ & $2$ & $0$ \\
Type $2$~\cite{park2014new}  & $1676$ & $3$ & $2$ & $2$ \\
Type $2$~\cite{RodriguezKoc} & $3430$  & $6$ & $4$ & $4$ \\
Type $3$~\cite{Masoleh:2004}  & $539$ & $0$ & $0$ & $0$ \\ \hline
Ours, $b \not = 2c$& $728$ & $2$ & $2$ & $1$ \\
Ours, $b=2c$ & $2$ & $0$ & $0$ & $0$ \\
\end{tabular}
\end{table}

\setlength{\tabcolsep}{1.5pt}
\renewcommand\arraystretch{0.65}
\setlength{\LTcapwidth}{\linewidth}
{\footnotesize
\begin{longtable}{l|l|l|l|l}
\caption{Our family of irreducible pentanomials and their number of XORs ($b,c,N_\oplus$), $2b\neq c$.}
\label{tab:penta:irredut:2}\\ \hline
$2,1,11$ & $3,2,22$ & $4,1,25$ & $5,1,31$ & $5,2,34$\\
$6,1,37$ & $5,3,37$ & $7,2,46$ & $9,5,67$ & $8,7,67$\\
$9,6,70$ & $12,1,73$ & $11,3,73$ & $10,7,79$ & $13,3,85$\\
$10,9,85$ & $13,4,88$ & $15,6,106$ & $14,9,109$ & $19,2,118$\\
$17,6,118$ & $15,10,118$ & $17,11,133$ & $17,12,136$ & $21,5,139$\\
$20,7,139$ & $16,15,139$ & $21,6,142$ & $23,5,151$ & $22,7,151$\\
$25,2,154$ & $21,11,157$ & $21,13,163$ & $27,5,175$ & $23,13,175$\\
$29,2,178$ & $25,10,178$ & $23,14,178$ & $25,12,184$ & $28,7,187$\\
$32,1,193$ & $28,9,193$ & $31,4,196$ & $23,20,196$ & $30,7,199$\\
$28,15,211$ & $27,18,214$ & $35,3,217$ & $31,11,217$ & $27,22,226$\\
$29,20,232$ & $35,10,238$ & $31,19,241$ & $38,7,247$ & $31,21,247$\\
$41,3,253$ & $38,9,253$ & $37,12,256$ & $35,19,265$ & $39,12,268$\\
$34,25,277$ & $45,4,280$ & $33,29,283$ & $47,2,286$ & $40,17,289$\\
$38,23,295$ & $48,7,307$ & $40,23,307$ & $46,15,319$ & $42,23,319$\\
$53,2,322$ & $45,18,322$ & $41,26,322$ & $45,19,325$ & $38,33,325$\\
$41,28,328$ & $52,7,331$ & $41,29,331$ & $47,20,340$ & $45,26,346$\\
$43,30,346$ & $49,19,349$ & $41,35,349$ & $45,28,352$ & $57,6,358$\\
$51,18,358$ & $45,30,358$ & $46,31,367$ & $55,14,370$ & $52,25,385$\\
$63,4,388$ & $62,7,391$ & $45,44,400$ & $51,34,406$ & $59,19,409$\\
$50,41,421$ & $63,18,430$ & $68,9,433$ & $63,19,433$ & $59,27,433$\\
$56,33,433$ & $67,12,436$ & $69,11,445$ & $60,31,451$ & $75,2,454$\\
$56,41,457$ & $63,29,463$ & $62,31,371$ & $59,37,463$ & $75,6,466$\\
$71,14,466$ & $65,26,466$ & $61,36,472$ & $77,5,475$ & $74,15,487$\\
$63,37,487$ & $67,30,490$ & $65,34,490$ & $73,19,493$ & $71,30,514$\\
$87,2,526$ & $87,6,538$ & $75,30,538$ & $69,42,538$ & $82,17,541$\\
$71,46,562$ & $70,49,565$ & $81,28,568$ & $77,36,568$ & $85,21,571$\\
$65,61,571$ & $83,28,580$ & $95,10,598$ & $85,30,598$ & $75,50,598$\\
$95,12,604$ & $98,9,613$ & $86,33,613$ & $81,43,613$ & $78,49,613$\\
$77,51,613$ & $103,3,625$ & $91,28,628$ & $87,37,631$ & $78,55,631$\\
$101,11,637$ & $74,65,637$ & $104,7,643$ & $81,54,646$ & $79,60,652$\\
$79,61,655$ & $101,18,658$ & $85,53,667$ & $112,1,673$ & $91,44,676$\\
$90,47,679$ & $79,69,679$ & $81,66,682$ & $105,19,685$ & $90,49,685$\\
$95,43,697$ & $79,75,697$ & $102,31,703$ & $99,37,703$ & $91,53,703$\\
$97,42,706$ & $94,49,709$ & $104,31,715$ & $119,2,718$ & $105,30,718$\\
$110,23,727$ & $103,37,727$ & $105,34,730$ & $99,46,730$ & $88,73,745$\\
$99,52,748$ & $118,15,751$ & $103,45,751$ & $95,61,751$ & $115,23,757$\\
$105,43,757$ & $93,67,757$ & $125,4,760$ & $93,68,760$ & $127,2,766$\\
$87,83,769$ & $123,14,778$ & $130,1,781$ & $97,67,781$ & $128,7,787$\\
$108,47,787$ & $103,59,793$ & $92,81,793$ & $119,30,802$ & $99,70,802$\\
$117,36,808$ & $120,31,811$ & $105,61,811$ & $119,34,814$ & $106,63,823$\\
$131,14,826$ & $133,13,835$ & $140,1,841$ & $95,91,841$ & $123,37,847$\\
$111,61,847$ & $115,54,850$ & $118,49,853$ & $113,59,853$ & $141,6,862$\\
$107,76,868$ & $130,31,871$ & $125,42,874$ & $125,43,877$ & $142,15,895$\\
$139,22,898$ & $125,50,898$ & $115,70,898$ & $131,43,913$ & $154,1,925$\\
$142,25,925$ & $155,3,937$ & $107,102,946$ & $154,9,949$ & $114,89,949$\\
$109,99,949$ & $145,34,970$ & $137,50,970$ & $135,54,970$ & $123,78,970$\\
$146,33,973$ & $145,36,976$ & $133,60,976$ & $121,85,979$ & $161,6,982$\\
$143,44,988$ & $123,84,988$ & $129,74,994$ & $153,29,1.003$ & $156,25,1009$\\
$115,107,1.009$ & $118,105,1.021$ & $169,4,1.024$ & $145,52,1.024$ & $137,68,1024$\\
$125,92,1.024$ & $139,67,1.033$ & $135,78,1.042$ & $129,90,1.042$ & $129,91,1045$\\
$135,84,1.060$ & $174,7,1.063$ & $126,103,1.063$ & $157,42,1.066$ & $161,35,1069$\\
$154,49,1.069$ & $133,93,1.075$ & $171,18,1.078$ & $153,54,1.078$ & $135,90,1078$\\
$179,5,1.087$ & $130,103,1.087$ & $169,27,1.093$ & $162,41,1.093$ & $142,81,1093$\\
$133,99,1.093$ & $122,121,1.093$ & $124,121,1.105$ & $130,113,1.117$ & $173,29,1123$\\
$167,43,1.129$ & $144,89,1.129$ & $189,4,1.144$ & $177,28,1.144$ & $161,60,1144$\\
$163,62,1.162$ & $133,123,1.165$ & $140,111,1.171$ & $147,101,1.183$ & $193,10,1186$\\
$185,27,1.189$ & $189,20,1.192$ & $197,6,1.198$ & $175,50,1.198$ & $160,81,1201$\\
$135,132,1.204$ & $170,63,1.207$ & $166,71,1.207$ & $149,109,1.219$ & $153,102,1222$\\
$191,28,1.228$ & $189,37,1.243$ & $161,93,1.243$ & $159,100,1.252$ & $179,61,1255$\\
$155,109,1.255$ & $203,14,1.258$ & $161,98,1.258$ & $198,25,1.261$ & $170,81,1261$\\
$150,121,1.261$ & $149,132,1.288$ & $205,21,1.291$ & $189,54,1.294$ & $163,109,1303$\\
$151,134,1.306$ & $173,93,1.315$ & $148,143,1.315$ & $209,22,1.318$ & $187,66,1318$\\
$196,49,1.321$ & $190,63,1.327$ & $183,77,1.327$ & $194,57,1.333$ & $172,105,1345$\\
$223,4,1.348$ & $173,108,1.360$ & $225,6,1.366$ & $204,49,1.369$ & $155,149,1375$\\
$162,137,1.381$ & $161,140,1.384$ & $204,55,1.387$ & $193,77,1.387$ & $199,69,1399$\\
$225,18,1.402$ & $213,42,1.402$ & $195,78,1.402$ & $197,76,1.408$ & $183,108,1420$\\
$234,7,1.423$ & $203,69,1.423$ & $209,59,1.429$ & $161,155,1.429$ & $235,10,1438$\\
$235,12,1.444$ & $179,124,1.444$ & $218,49,1.453$ & $169,147,1.453$ & $201,90,1474$\\
$225,44,1.480$ & $173,148,1.480$ & $220,63,1.507$ & $248,9,1.513$ & $247,12,1516$\\
$254,1,1.525$ & $213,90,1.546$ & $217,83,1.549$ & $201,115,1.549$ & $224,71,1555$\\
$238,47,1.567$ & $261,6,1.582$ & $183,163,1.585$ & $227,76,1.588$ & $218,95,1591$\\
$178,175,1.591$ & $265,4,1.600$ & $241,53,1.603$ & $196,143,1.603$ & $267,2,1606$\\
$269,2,1.618$ & $265,10,1.618$ & $261,18,1.618$ & $241,58,1.618$ & $225,90,1618$\\
$221,98,1.618$ & $207,126,1.618$ & $205,130,1.618$ & $246,49,1.621$ & $272,1,1633$\\
$196,153,1.633$ & $192,161,1.633$ & $203,140,1.636$ & $254,39,1.639$ & $194,161,1645$\\
$257,37,1.651$ & $212,127,1.651$ & $239,77,1.663$ & $255,46,1.666$ & $227,102,1666$\\
$245,67,1.669$ & $234,89,1.669$ & $197,163,1.669$ & $209,140,1.672$ & $244,71,1675$\\
$247,68,1.684$ & $195,172,1.684$ & $195,173,1.687$ & $213,138,1.690$ & $274,17,1693$\\
$193,180,1.696$ & $280,9,1.705$ & $215,139,1.705$ & $243,84,1.708$ & $218,135,1711$\\
$239,94,1.714$ & $219,134,1.714$ & $241,91,1.717$ & $216,145,1.729$ & $225,130,1738$\\
$223,134,1.738$ & $215,150,1.738$ & $249,84,1.744$ & $256,71,1.747$ & $208,167,1747$\\
$211,163,1.753$ & $231,124,1.756$ & $255,77,1.759$ & $199,189,1.759$ & $230,129,1765$\\
$213,163,1.765$ & $249,92,1.768$ & $295,2,1.774$ & $265,66,1.786$ & $255,86,1786$\\
$286,25,1.789$ & $285,30,1.798$ & $255,90,1.798$ & $225,150,1.798$ & $267,67,1801$\\
$263,75,1.801$ & $211,181,1.807$ & $293,18,1.810$ & $285,36,1.816$ & $247,116,1828$\\
$259,94,1.834$ & $266,81,1.837$ & $253,107,1.837$ & $221,171,1.837$ & $285,44,1840$\\
$300,17,1.849$ & $252,113,1.849$ & $279,61,1.855$ & $265,91,1.861$ & $249,124,1864$\\
$244,137,1.873$ & $227,172,1.876$ & $273,84,1.888$ & $252,127,1.891$ & $311,13,1903$\\
$271,93,1.903$ & $266,103,1.903$ & $259,117,1.903$ & $265,109,1.915$ & $255,131,1921$\\
$252,137,1.921$ & $215,212,1.924$ & $298,47,1.927$ & $231,181,1.927$ & $305,36,1936$\\
$245,157,1.939$ & $323,2,1.942$ & $243,162,1.942$ & $259,131,1.945$ & $223,203,1945$\\
$279,92,1.948$ & $238,175,1.951$ & $274,105,1.957$ & $325,6,1.966$ & $292,73,1969$\\
$322,15,1.975$ & $319,22,1.978$ & $303,54,1.978$ & $253,154,1.978$ & $310,47,1999$\\
$329,14,2.014$ & $314,47,2.023$ & $323,30,2.026$ & $257,162,2.026$ & $314,49,2029$\\
$323,34,2.038$ & $289,102,2.038$ & $255,170,2.038$ & $307,68,2.044$ & $243,198,2050$\\
$329,27,2.053$ & $253,179,2.053$ & $237,211,2.053$ & $256,175,2.059$ & $339,11,2065$\\
$308,73,2.065$ & $303,83,2.065$ & $243,203,2.065$ & $287,116,2.068$ & $243,205,2071$\\
$266,161,2.077$ & $305,91,2.101$ & $320,63,2.107$ & $301,101,2.107$ & $343,19,2113$\\
$243,220,2.116$ & $293,122,2.122$ & $349,11,2.125$ & $285,139,2.125$ & $253,203,2125$\\
$266,183,2.143$ & $254,207,2.143$ & $307,102,2.146$ & $325,69,2.155$ & $357,6,2158$\\
$315,90,2.158$ & $349,26,2.170$ & $329,67,2.173$ & $340,49,2.185$ & $347,37,2191$\\
$341,50,2.194$ & $297,138,2.194$ & $285,164,2.200$ & $283,173,2.215$ & $270,199,2215$\\
$349,42,2.218$ & $301,139,2.221$ & $301,141,2.227$ & $261,221,2.227$ & $365,18,2242$\\
$297,156,2.248$ & $365,21,2.251$ & $268,217,2.257$ & $371,13,2.263$ & $371,14,2266$\\
$287,182,2.266$ & $374,9,2.269$ & $361,36,2.272$ & $328,103,2.275$ & $375,10,2278$\\
$260,241,2.281$ & $279,204,2.284$ & $313,139,2.293$ & $257,251,2.293$ & $297,173,2299$\\
$264,239,2.299$ & $381,6,2.302$ & $304,161,2.305$ & $260,249,2.305$ & $355,62,2314$\\
$321,130,2.314$ & $372,31,2.323$ & $341,93,2.323$ & $293,189,2.323$ & $364,49,2329$\\
$287,203,2.329$ & $351,76,2.332$ & $377,26,2.338$ & $369,42,2.338$ & $325,130,2338$\\
$299,182,2.338$ & $378,25,2.341$ & $321,140,2.344$ & $347,91,2.353$ & $332,121,2353$\\
$361,66,2.362$ & $303,182,2.362$ & $278,233,2.365$ & $305,187,2.389$ & $392,15,2395$\\
$311,180,2.404$ & $386,31,2.407$ & $271,261,2.407$ & $395,14,2.410$ & $307,190,2410$\\
$297,210,2.410$ & $320,169,2.425$ & $351,108,2.428$ & $389,35,2.437$ & $361,93,2443$\\
$357,102,2.446$ & $404,9,2.449$ & $343,133,2.455$ & $287,245,2.455$ & $403,14,2458$\\
$335,150,2.458$ & $325,170,2.458$ & $293,234,2.458$ & $397,27,2.461$ & $286,255,2479$\\
$393,42,2.482$ & $365,101,2.491$ & $395,44,2.500$ & $411,14,2.506$ & $283,270,2506$\\
$381,76,2.512$ & $397,45,2.515$ & $285,269,2.515$ & $321,203,2.533$ & $407,38,2554$\\
$299,254,2.554$ & $321,211,2.557$ & $336,185,2.569$ & $320,217,2.569$ & $411,38,2578$\\
$403,54,2.578$ & $355,150,2.578$ & $339,182,2.578$ & $322,217,2.581$ & $423,18,2590$\\
$403,59,2.593$ & $389,91,2.605$ & $358,153,2.605$ & $321,228,2.608$ & $320,231,2611$\\
$379,115,2.617$ & $425,27,2.629$ & $389,99,2.629$ & $353,173,2.635$ & $435,10,2638$\\
$400,81,2.641$ & $396,89,2.641$ & $351,181,2.647$ & $326,231,2.647$ & $295,294,2650$\\
$422,41,2.653$ & $382,121,2.653$ & $363,164,2.668$ & $319,252,2.668$ & $303,284,2668$\\
$311,270,2.674$ & $401,91,2.677$ & $325,243,2.677$ & $373,148,2.680$ & $443,14,2698$\\
$417,66,2.698$ & $413,74,2.698$ & $375,150,2.698$ & $345,210,2.698$ & $301,298,2698$\\
$362,177,2.701$ & $381,140,2.704$ & $364,175,2.707$ & $443,19,2.713$ & $367,173,2719$\\
$405,98,2.722$ & $448,17,2.737$ & $375,163,2.737$ & $407,102,2.746$ & $405,106,2746$\\
$377,162,2.746$ & $427,67,2.761$ & $316,289,2.761$ & $439,45,2.767$ & $339,245,2767$\\
$318,287,2.767$ & $461,4,2.776$ & $393,140,2.776$ & $457,13,2.779$ & $445,37,2779$\\
$423,83,2.785$ & $403,124,2.788$ & $335,262,2.794$ & $413,107,2.797$ & $392,151,2803$\\
$344,249,2.809$ & $387,166,2.818$ & $355,230,2.818$ & $389,164,2.824$ & $466,15,2839$\\
$362,223,2.839$ & $321,306,2.842$ & $353,243,2.845$ & $462,31,2.863$ & $411,133,2863$\\
$394,169,2.869$ & $441,76,2.872$ & $436,89,2.881$ & $338,287,2.887$ & $443,78,2890$\\
$373,218,2.890$ & $421,123,2.893$ & $480,7,2.899$ & $380,207,2.899$ & $435,102,2914$\\
$411,150,2.914$ & $405,162,2.914$ & $369,234,2.914$ & $376,223,2.923$ & $420,137,2929$\\
$435,108,2.932$ & $399,180,2.932$ & $458,63,2.935$ & $445,89,2.935$ & $354,271,2935$\\
$437,107,2.941$ & $401,179,2.941$ & $425,133,2.947$ & $483,18,2.950$ & $350,287,2959$\\
$429,132,2.968$ & $369,252,2.968$ & $397,197,2.971$ & $392,207,2.971$ & $364,265,2977$\\
$494,7,2.983$ & $387,222,2.986$ & $494,9,2.989$ & $429,139,2.989$ & $475,50,2998$\\
$425,150,2.998$ & $375,250,2.998$ & $431,140,3.004$ & $466,71,3.007$ & $419,165,3007$\\
$337,332,3.016$ & $427,156,3.028$ & $407,196,3.028$ & $347,316,3.028$ & $487,37,3031$\\
$457,98,3.034$ & $355,302,3.034$ & $485,43,3.037$ & $365,284,3.040$ & $415,187,3049$\\
$418,183,3.055$ &  & \\\hline

\end{longtable}
}

\section{Conclusions} \label{conclusions}
In this paper, we present a new class of pentanomials over $\mathbb{F}_2$, defined
by $x^{2b+c}+ x^{b+c}+ x^b + x^c + 1$. We give the exact number
of XORs in the reduction process; we note that in the reduction process
no ANDs are required.

It is interesting to point out that even though the cases $c=1$ and $c>1$,
as shown in Figures~\ref{fig:Fig-DRed-ka2}~and~\ref{fig:Fig-DRed-ka3}, are
quite different, the final result in terms of number of XORs is the same.
We also consider a special case when $b=2c$ where further reductions are
possible.

There are irreducible pentanomials of this shape for several degree extensions
of practical interest. We provide a detailed analysis of the space and time
complexity involved in the reduction using the pentanomials in our family.
For the multiplication process, we simply use the standard Karatsuba algorithm.

The proved complexity analysis of the multiplier and reduction considering
the family proposed in this paper, as well as our analysis
suggests that these pentanomials are as good as or possibly better to the
ones already proposed.

We leave for future work to produce a one-step algorithm using
our pentanomials, that is, a multiplier that performs multiplication
and reduction in a single step using our family of polynomials, as
well as a detailed study of the delay obtained using this algorithm.

\bibliographystyle{plain}

\end{document}